\documentclass[12pt]{article}
\usepackage{amsthm,amsmath,amssymb,amscd,amsfonts,mathrsfs,color,graphicx,bm,framed}
\begin{document}

\medbreak\centerline{\bf The Role of the Gradient Term of the Bochner-Kodaira Formula}
\centerline{\bf in Coherent Sheaf Extension}

\medbreak
\centerline{\sc Yum-Tong Siu}

\medbreak
\centerline{\it In Memory of Nessim Sibony (1947 - 2021)}

\vskip.35in
{\footnotesize\bigbreak\noindent{\sc Abstract.}  In applying the Bochner-Kodaira formula with boundary term to solve the $\bar\partial$ equation with $L^2$ estimates, the gradient term is usually not used.  Two potentially important applications of the use of the gradient term are the strong rigidity for holomorphic vector bundles and the very ampleness part of the Fujita conjecture. In this note we use the gradient term to construct holomorphic sections to prove the Thullen-type extension across codimension $1$ for holomorphic vector bundles with Hermitian metric whose curvature is $L^p$ for some $p>1$.  This construction of sections points out a typical way of how the gradient term can be used.

}
\vskip.25in
\noindent{\bf\S1}\ \  {\sc Introduction}

\bigbreak A very powerful tool to produce holomorphic sections of holomorphic vector bundles, first introduced by Kodaira, is to use the positivity of the curvature tensor in the Bochner-Kodaira formula.  The tool, with additional techniques from Morrey, Kohn and H\"ormander, can be applied to the case of pseudoconvex boundary.  Most of the time the gradient term in the formula is not used except for its nonnegativity.  In the case of nonnegative curvature, the gradient term can be used to yield a foliation and to construct Kohn's multiplier ideal sheaves [Koh].  In some settings, even when the curvature tensor changes sign, the gradient term can be used for the construction of holomorphic sections.  Two potentially important applications of the use of the gradient term are the strong rigidity for holomorphic vector bundles and the very ampleness part of the Fujita conjecture.  Both will be explained in more details below in (1.4).
In this note we use the gradient term to construct holomorphic sections to get the following Thullen-type extension across codimension $1$ for holomorphic vector bundles whose curvature is $L^p$ for some $p>1$.  This construction of sections points out a typical way of how the gradient term can be used.  To use the gradient term which is for the $(0,1)$ directions, the Cauchy integral formula for smooth functions on local complex curves and the estimate for singular integrals, especially potential estimates, are used in (2.2) and (2.1) below.

\bigbreak\noindent(1.1) {\it Main Theorem.}  Let $X$ be complex manifold of complex dimension $n\geq 2$ and $Y$ be a connected nonsingular complex hypersurface in $X$ and $G$ be an open subset of $X$ which intersects $Y$.  Let $E$ be a holomorphic vector bundle on $(X-Y)\cup G$ and $h$ be a smooth Hermitian metric of $E$ such that the pointwise norm of the curvature $\Theta_h$ of $h$, with respect to $h$ and any smooth Hermitian metric of $X$, is locally $L^p$ for some $p>1$ as a function on $X$.  Then $E$ can be extended to a coherent sheaf ${\mathcal F}$ on $X$, which is reflexive and therefore unique.  Moreover, ${\mathcal F}$ can be defined by the presheaf which assigns to an open subset $Q$ of $X$ the vector space ${\mathcal P}_Q$ with the following definition. An element $s$ of ${\mathcal P}_Q$ is a holomorphic section of $Q-Y$ such that for some (possibly empty) complex subvariety $Z_s$ of complex dimension $\leq n-2$ in $Q\cap Y$, the pointwise norm $|s|$ of $s$ as a function of $Q-Y$ is locally $L^2$ at points of $Q-Y-Z_s$ as a function of $Q$.

\bigbreak The special case of extending a holomorphic vector bundle with $L^2$ curvature tensor across complex subvarieties of codimension $\geq 2$ was already done in [Ban-Siu], by putting together Bando's result of extension across a single point [Ban] and the slicing result in Hartogs extension of coherent sheaves [Siu3].  Besides establishing the extension with $L^p$ curvature for $p>1$ in the natural more general setting, this note presents a more straightforward and more transparent use of the Bochner-Kodaira formula with boundary term, instead of using the solution of the {\it Dirichlet-$\bar\partial$-Neumann problem} as in the arguments in [Ban].  It is precisely this straightforward use of the Bochner-Kodaira formula with boundary term that makes it possible to provide the more general result in the Main Theorem.

\bigbreak For $r>0$ let ${\mathbb D}_r$ be the open disk in ${\mathbb C}$ centered at the origin with radius $r$. Denote ${\mathbb D}_1$ by ${\mathbb D}$.  We also use the notations ${\mathbb D}_r^*={\mathbb D}_r-\{0\}$ and ${\mathbb D}^*={\mathbb D}-\{0\}$.  The notation ${\mathcal O}$ denotes the sheaf of germs of holomorphic functions and ${\mathcal O}(E)$ denotes the locally free sheaf associated to the holomorphic vector bundle $E$.

\bigbreak\noindent(1.2) {\it Perspective of Extension of Analytic Objects.}  To see how the extension result of the Main Theorem fits into the general question of extension of analytic objects and to indicate the possible directions of further development, we enumerate here the more common settings for extension, analytic objects to be extended and conditions on them.

\bigbreak\noindent Here are the more common settings for extension (where the use of polydisks can alternatively be replaced by balls).

\begin{itemize}
\item[(1)] Extension from the Hartogs's figure $\left(({\mathbb D}^k-\overline{{\mathbb D}_a^k})\times{\mathbb D}^{n-k}\right)\cup({\mathbb D}^k\times{\mathbb D}_b^{n-k})$ of codimension $k$ to ${\mathbb D}^n$ for $a,b\in(0,1)$ and $n>k$.
\item[(2)] The Thullen setting from $\left(({\mathbb D}^k-\{0\})\times{\mathbb D}^{n-k}\right)\cup({\mathbb D}^k\times{\mathbb D}_b^{n-k})$ of codimension $k$ to ${\mathbb D}^n$ for $0<b<1$ and $n>k$.  This includes extension across subvarieties of codimension $\geq k+1$ as a consequence.
\item[(3)] Extension across submanifold of codimension $k$ from $({\mathbb D}^k-\{0\})\times{\mathbb D}^{n-k}$ to ${\mathbb D}^n$.
\end{itemize}

\noindent Here are the more common analytic objects to be extended.
\begin{itemize}
\item[(i)] Holomorphic functions.
\item[(ii)] Subvarieties.
\item[(iii)] Closed positive currents.
\item[(iv)] Holomorphic line or vector bundles.
\item[(v)] Coherent analytic sheaves or subsheaves.
\item[(vi)] Holomorphic or meromorphic maps to complex Hermitian or K\"ahler manifolds.
\end{itemize}

\noindent Here are the more common conditions on the analytic objects to be extended.
\begin{itemize}
\item[(a)] No additional condition other than holomorphicity.
\item[(b)] Nonnegativity condition
\item[(c)] Condition of finite mass or $L^p$ norm.
\end{itemize}

\noindent For example, holomorphic functions without any additional assumption have Hartogs extension of codimension $\geq 1$.  Holomorphic functions which are $L^2$ can be extended across a complex submanifold of codimension $\geq 1$.  Subvarieties of pure complex codimension $k$ have Hartogs extension of codimension $\geq k$.  Closed positive $(1,1)$-currents of finite mass can be extended across complex submanifold of codimension $\geq 1$.  Coherent analytic sheaves, with the property that local holomorphic sections extend uniquely across local subvarieties of complex codimension $\geq k$, have Hartogs extension of codimension $\geq k$.  Holomorphic vector bundles, as special cases of coherent analytic sheaves, have Hartogs extension of codimension $\geq 2$.  Details about such and other results of extension of analytic objects can be found in
[Ban, Ban-Siu, Bis, Harv, Hart, Shi1, Shi2, Sib, Siu2, Siu3, Siu4, Siu-Tra, Sko, Thu, Tra1, Tra2]. The Main Theorem in this note deals with the case of Thullen extension of codimension $\geq 1$ under the additional condition on the $L^p$ integrability of the curvature tensor of a Hermitian metric for some $p>1$.  We would like to remark that besides the settings for extension enumerated above which are formulated in terms of subvarieties and polydisks or balls, there are also more special settings for extension, such as across pluripolar sets, across sets with conditions on their Hausdorff measure, across totally real sets or across sets defined by more general conditions involving CR manifolds (see for example, the extension of closed positive currents for some such special settings in [Sib]).

\bigbreak\noindent(1.3) {\it Comments on Directions of Further Development.}  The question naturally arises about the condition corresponding to $L^p$ integrability of curvature for $p>1$ when the holomorphic vector bundle $E$ in the Main Theorem is replaced by a coherent analytic sheaf.  This question belongs to the second one of the two natural directions of further development along the lines of the Main Theorem we would like to comment on.

\bigbreak One is the weakening of the integrability condition for curvature in the Main Theorem.  Condition of integrability of $|\Theta_h|^p$ for some $p>1$ comes from the integrability of $$\int_{z\in{\mathbb C},\,|z|<1}\frac{1}{|z|^{2-\varepsilon}}$$ for some $\varepsilon>0$.  A weakening of the condition to the integrability of the product of $|\Theta_h|(\log|\Theta_h|)^q$ for an appropriate positive $q$ is expected from the integrability of
$$\int_{z\in{\mathbb C},\,|z|<1}\frac{1}{|z|^2(\log|z|^2)^2}.$$

\bigbreak Another is to regard $E$ over $X-Y$ as a holomorphic family of vector spaces which can be generalized to more general holomorphic families, for example, of linear fiber spaces associated to coherent analytic sheaves as defined in [Fis, Gra, Gro], or even more general complex spaces.  The development will be in the setting of extending a holomorphic family of complex spaces from over the parameter space $(X-Y)\cup G$ to over the parameter $X$ with integrability condition on tensors from appropriately defined metrics.  These two natural directions of further development will not be discussed here.

\bigbreak\noindent(1.4) {\it Two Potentially Important Applications of the Use of the Gradient term.}  One of the main statements in the strong rigidity theory for compact K\"ahler manifolds [Siu1] is that if $X$ is a compact K\"ahler manifold of complex dimension $\geq 2$ whose curvature satisfies some appropriate curvature condition, then any compact K\"ahler manifold $Y$ which is of the same homotopy type as $X$ must be biholomorphic or conjugate biholomorphic to $X$.  The key argument is to apply a form of the Bochner-Kodaira formula to a harmonic map from $Y$ to $X$ which is a homotopy equivalence.  A natural question is to ask for an analogous strong rigidity theory for holomorphic vector bundles $E$ over a fixed compact K\"ahler (or even algebraic) manifold $M$.

\bigbreak One seeks an appropriate condition on $E$ so that a holomorphic vector bundle $F$ over $M$ which is topologically equivalent to $E$ should be biholomorphic or conjugate holomorphic to $E$.  The obstacle is that a condition for the infinitesimal rigidity of $E$ is the vanishing of $H^1(M,E\otimes E^*)$ which cannot be obtained as a consequence of the curvature term in the Bochner-Kodaira formula, because the curvature is for the tensor product of $E$ and its dual $E^*$.  For such an approach, the strong rigidity theory or even the infinitesimal rigidity for $E$ must be from the use of the gradient term of the Bochner-Kodaira formula.  How it is to be done remains an open problem.

\bigbreak The second potentially important application of the use the gradient term is the very ampleness part of the Fujita conjecture [Fuj].  For a compact algebraic manifold $X$ of complex dimension $n$ and an ample line bundle $L$ over $X$, the original Fujita conjecture seeks to prove that $mL+K_X$ is globally free for $m\geq n+1$ and is very ample for $m\geq n+2$.  The freeness part is known under the stronger assumption of $m\geq\frac{1}{2}n(n+1)+1$ [Ang-Siu] which can be improved to $m\geq\left(e+\frac{1}{2}\right)n^{\frac{4}{3}}+\frac{1}{2}n^{\frac{2}{3}}+1$ [Hei] but not to the original conjecture bound of $m\geq n+1$.  For the very ampleness part, the very ampleness of $mL+2K_X$ for $m\geq m_n$ is known for some positive integer $m_n$ which depends only on $n$ [Dem1, Dem2, Siu5].  Thus for any positive rational number $\varepsilon$ and any prescribed jets at a finite set $Z$ of prescribed points in $X$ with total jet order $\leq q$, there exists some positive integer $m_{n,q,\varepsilon}$ which depends only on $n, q, \varepsilon$ with the property that for $m\geq m_{n,q,\varepsilon}$ there exists a metric $h$ for the ${\mathbb Q}$-bundle $mL+\varepsilon K_X$ over $X$ which is smooth outside $Z$ whose curvature current is positive such that the pole-order of the metric at a point $P$ of $Z$ is at least the prescribed jet order at $P$ plus $n$.  If the metric $h$ is for $mL$ instead of for $mL+\varepsilon K_X$, then by Nadel's vanishing theorem [Nad] one can construct from the metric a holomorphic section of $mL+K_X$ which assumes the prescribed jets at points of $Z$.  Now, with the metric $h$ for $mL+\varepsilon K_X$, for Nadel's vanishing theorem one needs to use the gradient term in the Bochner-Kodaira formula to handle $\varepsilon$ times the Ricci curvature in such a way that the argument depends only on the complex dimension of $X$ and not on $X$ itself.  How it is to be done remains an open problem.

\bigbreak The rest of this note is devoted to proving the Main Theorem, which will be obtained as a corollary to a slightly more general Theorem(5.5).  The details of the proof seem a bit involved, because all the steps and the constants used are explicitly presented to make it more convenient for the reader to follow.  Actually, the basic argument is just a simple use of the Cauchy integral formula for smooth functions in single complex variable, together with H\"older-type inequalities and integration by parts, to absorb the curvature term into the gradient term in the Bochner-Kodaira formula with boundary term.

\vskip.2in
\bigbreak\noindent{\bf\S2}\ \  {\sc Potential Estimate for Cauchy Kernel}

\bigbreak\noindent(2.1) {\it Potential Estimate of Singular Kernel.}  Let $\Omega$ be a bounded domain in ${\mathbb R}^m$ (with coordinate $x$) and $h(x,y)$ be a nonnegative measurable function on $\Omega\times\Omega$.  For $p\geq 1$ define
$$
(V_h f)(x)=\int_{y\in\Omega}h(x,y)f(y)dy
$$
for any $L^p$ function $f$ on $\Omega$.
Let $p\leq q\leq\infty$ and
$$
\frac{1}{r}=1+\frac{1}{q}-\frac{1}{p}.
$$
Assume that the $L^r(\Omega)$ norm of $h(x-y)$ as a function of $y$ is uniformly bounded for fixed $x$ in $\Omega$.
Then $V_h$ maps $L^p(\Omega)$ to $L^q(\Omega)$
and
$$
\|V_h f\|_{L^q(\Omega)}\leq\left(\sup_{x\in\Omega}\|h(x,\cdot)\|_{L^r(\Omega)}\right)\|f\|_{L^p(\Omega)}.
$$

\bigbreak This statement is is just a trivial straightforward reformulation of [Gil-Tru, Lemma 7.12].  Let $\delta=\frac{1}{p}-\frac{1}{q}$.
The following verification is simply a reproduction of the argument of the proof of [Gil-Tru, Lemma 7.12], which applies H\"older's inequality, for three factors with respective exponents $$q,\quad\frac{1}{1-\frac{1}{p}},\quad{\rm and}\quad\frac{1}{\delta}$$ whose reciprocals add up to $1$ in
$$
\frac{1}{q}+\left(1-\frac{1}{p}\right)+\delta=1,
$$
to
$$
h(x,y)|f(y)|=\left(h(x,y)^r|f(y)|^p\right)^{\frac{1}{q}}\left(h(x,y)^r\right)^{1-\frac{1}{p}}\left(|f(y)|^p\right)^\delta,
$$
where among the exponents on the right-hand side
$$
\frac{r}{q}+r\left(1-\frac{1}{p}\right)=1\quad{\rm and}\quad\frac{p}{q}+p\,\delta=1
$$
hold.
The application of H\"older's inequality yields
$$
\begin{aligned}\left|(V_h f)(x)\right|&\leq\int_{y\in\Omega}h(x,y)\left|f(y)\right|dy\cr
&=\int_{y\in\Omega}\left(h(x,y)^r|f(y)|^p\right)^{\frac{1}{q}}\left(h(x,y)^r\right)^{1-\frac{1}{p}}\left(|f(y)|^p\right)^\delta dy\cr
&\leq\left(\int_{y\in\Omega}h(x,y)^r|f(y)|^p dy\right)^{\frac{1}{q}}\left(\int_{y\in\Omega}h(x,y)^r dy\right)^{1-\frac{1}{p}}
\left(\int_{y\in\Omega}
|f(y)|^p dy\right)^\delta\cr
\end{aligned}
$$
so that if we let
$$
A=\sup_{x\in\Omega}\int_{y\in\Omega}h(x,y)^r dy,
$$
then
$$
\begin{aligned}\left|(V_h f)(x)\right|^q&\leq\left(\int_{y\in\Omega}h(x,y)^r|f(y)|^p dy\right)\left(\int_{y\in\Omega}h(x,y)^r dy\right)^{q\left(1-\frac{1}{p}\right)}
\left(\int_{y\in\Omega}
|f(y)|^p dy\right)^{q\delta}\cr
&\leq\left(\int_{y\in\Omega}h(x,y)^r|f(y)|^p dy\right)A^{q\left(1-\frac{1}{p}\right)}
\left(\int_{y\in\Omega}
|f(y)|^p dy\right)^{q\delta}\cr
\end{aligned}
$$
and
$$
\begin{aligned}&\int_{x\in\Omega}\left|(V_\mu f)(x)\right|^q dx
\leq A^{q\left(1-\frac{1}{p}\right)}\left(\int_{x,y\in\Omega}h(x,y)^r|f(y)|^p dxdy\right)
\left(\int_{y\in\Omega}
|f(y)|^p dy\right)^{q\delta}\cr
&\qquad= A^{q\left(1-\frac{1}{p}\right)}\left(\int_{y\in\Omega}\left(\int_{x\in\Omega}h(x,y)^r dx\right)|f(y)|^p dy\right)
\left(\int_{y\in\Omega}
|f(y)|^p dy\right)^{q\delta}\cr
&\qquad\leq A^{q\left(1-\frac{1}{p}\right)}\left(\sup_{y\in\Omega}\left(\int_{x\in\Omega}h(x,y)^r dx\right)\right)\left(\int_{y\in\Omega}|f(y)|^p dy\right)
\left(\int_{y\in\Omega}
|f(y)|^p dy\right)^{q\delta}\cr
&\qquad= A^{1+q\left(1-\frac{1}{p}\right)}
\left(\int_{y\in\Omega}
|f(y)|^p dy\right)^{1+q\delta}\cr
\end{aligned}
$$
and
$$
\begin{aligned}\left\|V_h f\right\|_{L^q(\Omega)}&=\left(\int_{x\in\Omega}\left|(V_h f)(x)\right|^q dx\right)^{\frac{1}{q}}\cr
&\leq A^{\frac{1}{q}+\left(1-\frac{1}{p}\right)}
\left(\int_{y\in\Omega}
|f(y)|^p dy\right)^{\frac{1}{q}+\delta}\cr
&= A^{\frac{1}{r}}\left\|f\right\|_{L^p(\Omega)}.\cr
\end{aligned}
$$

\vskip.2in
\bigbreak For us the potential estimate will be applied to ${\mathbb C}^n$ instead of ${\mathbb R}^m$ and to the singular function
$$
h(z,\zeta)=\frac{1}{\pi|z_1-\zeta_1|},
$$
where $z=(z_1,\cdots,z_n),\,\zeta=(\zeta_1,\cdots,\zeta_n)$ are points in the bounded domain $\Omega$ of ${\mathbb C}^n$.
The use of the singular function $h(z,\zeta)$ comes from the Cauchy kernel for $(0,1)$-derivative with respect to the first coordinate $z_1$ of of ${\mathbb C}^n$ in the following Cauchy integral formula for (nonholomorphic) smooth functions of a single complex variable.

\vskip.2in
\bigbreak\noindent\noindent(2.2) {\it Cauchy's Integral Formula for Smooth Functions.}  Let $D$ be a smooth bounded domain in ${\mathbb C}$  with coordinate $z_1$.  Cauchy's integral formula for (nonholomorphic) smooth functions applied to a (nonholomorphic) function $\Phi(z_1)$ on $\bar D$ is
$$
\Phi(z_1)=\frac{1}{2\pi i}\int_{\zeta_1\in\partial D}\frac{\Phi(\zeta_1)d\zeta_1}{\zeta_1-z_1}+\frac{1}{2\pi i}\int_{\zeta_1\in D}\frac{(\partial_{\bar\zeta_1}\Phi(\zeta_1))d\zeta_1\wedge d\bar\zeta_1}{\zeta_1-z_1}
$$
for $z_1\in D$.  Suppose $\Phi(z_1)$ vanishes on $\partial D$.  Then the boundary term vanishes and
$$
\Phi(z_1)=\frac{1}{2\pi i}\int_{\zeta_1\in D}\frac{(\partial_{\bar\zeta_1}\Phi(\zeta_1))d\zeta_1\wedge d\bar\zeta_1}{\zeta_1-z_1}
$$
for $z_1\in\partial D$, which implies that
$$
\left|\Phi(z_1)\right|\leq\frac{1}{\pi}\int_{\zeta_1\in D}\frac{1}{|\zeta_1-z_1|}\left|\partial_{\bar\zeta_1}\Phi(\zeta_1)\right|
$$
for $z_1\in D$.

\vskip.3in
\bigbreak\noindent{\bf\S3}\ \  {\sc Relation of Two Types of Gradient Terms and Iteration of Potential Estimate}

\bigbreak\noindent(3.1) {\it Integration by Parts Applied to Covariant Differentiation of Sections of Holomorphic Vector Bundle with Hermitian Metric.}  Let $D$ be a bounded domain in ${\mathbb C}$ with smooth boundary and $W$ be a bounded domain in ${\mathbb C}^{n-1}$.  Let $E$ be a holomorphic vector bundle with a smooth Hermitian metric $h=\left(h_{\alpha\bar\beta}\right)$ on some open neighborhood of the topological closure $\bar D\times\bar W$ of $D\times W$ in ${\mathbb C}^n$ (with coordinates $z_1,\cdots,z_n$).  Let $f=(f^\alpha)$ be a smooth section of $E$ on $\bar D\times\bar W$ which vanishes on $(\partial D)\times\bar W$, where the superscript $\alpha$ denotes the index for the fiber coordinate of $E$.  Denote by $\nabla_{\bar z_1}$ and $\nabla_{z_1}$ the covariant differentiations of sections of $E$ with respect to $\bar z_1$ and $z_2$ respectively.  Though $\nabla_{z_1}$ depends on the connection of $h$, yet $\nabla_{\bar z_1}$ is simply $\bar\partial_{z_1}$, which is independent of the connection of $h$.
Let
$$
\omega=\left(\frac{i}{2}\,dz_2\wedge d\bar z_2\right)\wedge\cdots\wedge\left(\frac{i}{2}\,dz_n\wedge d\bar z_n\right).
$$
The vanishing of $f$ on $(\partial D)\times\bar W$ makes it possible to have the following integration by parts without boundary terms.
$$
\begin{aligned}&
\left\|\nabla_{\bar z_1}f\right\|_{D\times W}^2=\int_{D\times W}\sum_{\alpha,\beta} h_{\alpha\bar\beta}\left(\partial_{\bar z_1}f^\alpha\right)\overline{\left(\partial_{\bar z_1}f^\beta\right)}\left(\frac{i}{2}\,dz_1\wedge d\bar z_1\right)\wedge\omega\cr
&=\int_{(\partial D)\times W}\sum_{\alpha,\beta} h_{\alpha\bar\beta}\left(\partial_{\bar z_1}f^\alpha\right)\overline{f^\beta}\left(\frac{i}{2}\,d\bar z_1\right)\wedge\omega\cr
&-\int_{D\times W}\sum_{\alpha,\beta} h_{\alpha\bar\beta}\left(\nabla_{z_1}\bar\nabla_{z_1}f^\alpha\right)\overline{f^\beta}\left(\frac{i}{2}\,dz_1\wedge d\bar z_1\right)\wedge\omega\cr
&=
-\int_{D\times W}\sum_{\alpha,\beta} h_{\alpha\bar\beta}\left(\left[\nabla_{z_1},\,\bar\nabla_{z_1}\right]f^\alpha\right)\overline{f^\beta}\left(\frac{i}{2}\,dz_1\wedge d\bar z_1\right)\wedge\omega\cr
&-\int_{D\times W}\sum_{\alpha,\beta} h_{\alpha\bar\beta}\left(\bar\nabla_{z_1}\nabla_{z_1}f^\alpha\right)\overline{f^\beta}\left(\frac{i}{2}\,dz_1\wedge d\bar z_1\right)\wedge\omega\cr
&=-\int_{D\times W}\sum_{\alpha,\beta} \Theta_{\alpha\bar\beta 1\bar 1}\,f^\alpha\,\overline{f^\beta}\left(\frac{i}{2}\,dz_1\wedge d\bar z_1\right)\wedge\omega\cr
&+\int_{{\partial D}\times W}\sum_{\alpha,\beta} h_{\alpha\bar\beta}\left(\nabla_{z_1}f^\alpha\right)\overline{f^\beta}\left(\frac{i}{2}\,dz_1\right)\wedge\omega\cr
&+\int_{D\times W}\sum_{\alpha,\beta} h_{\alpha\bar\beta}\left(\nabla_{z_1}f^\alpha\right)\overline{\left(\nabla_{z_1}f^\beta\right)}\left(\frac{i}{2}\,dz_1\wedge d\bar z_1\right)\wedge\omega\cr
&=-\left(\Theta_{1\bar 1}f,\,f\right)_{D\times W}+\left\|\nabla_{z_1}f\right\|_{D\times W}^2.\cr
\end{aligned}
$$
Here $\Theta=\left(\Theta_{\alpha\bar\beta j\bar k}\right)$ is the curvature tensor for $h$ and $\Theta_{1\bar 1}=\left(\Theta_{\alpha\bar\beta 1\bar 1}\right)$ is its value at $\partial_{z_1}$ and $\partial_{\bar z_1}$ and defines a Hermitian form on the fibers of $E$.  For notational simplicity we have dropped the subscript $h$ of $\Theta_h$ and use $\Theta$ for the curvature tensor of $h$.

\vskip.2in
\bigbreak\noindent\bigbreak\noindent(3.2) {\it Iterated Application of Potential Estimate.}
Let $\gamma\geq 1$.  We are going to apply the potential estimate to the $\partial_{\bar z_1}$ of $|f|^2=\sum_{\alpha,\beta}h_{\alpha\bar\beta}f^\alpha\overline{f^\beta}$.
From
$$
\begin{aligned}\partial_{\bar z_1}\left(|f|^2\right)^{\gamma}&=\gamma\left(|f|^2\right)^{\gamma-1}\partial_{\bar z_1}|f|^2\cr
&=\gamma\left(|f|^2\right)^{\gamma-1}\left(\left<\nabla_{\bar z_1}f,\,f\right>+\left<f,\,\nabla_{z_1}f\right>\right)\cr
&=\gamma\,\left<\nabla_{\bar z_1}f,\,\left(|f|^2\right)^{\gamma-1} f\right>+\gamma\,\left<\left(|f|^2\right)^{\gamma-1} f,\,\nabla_{z_1}f\right>\cr
\end{aligned}
$$
it follows from H\"older's inequality that
$$
\begin{aligned}&\int_{D\times W}\left|\partial_{\bar z_1}\left(|f|^2\right)^{\gamma}\right|\leq \gamma\int_{D\times W}\left|\left<\nabla_{\bar z_1}f,\,\left(|f|^2\right)^{\gamma-1} f\right>\right|+\gamma\int_{D\times W}\left|\left<\left(|f|^2\right)^{\gamma-1} f,\,\nabla_{z_1}f\right>\right|\cr
&\leq \gamma\left\|\nabla_{\bar z_1}f\right\|_{L^2(D\times W)}\,\left\|\left(|f|^2\right)^{\gamma-1} f\right\|_{L^2(D\times W)}+\gamma\left\|\left(|f|^2\right)^{\gamma-1} f\right\|_{L^2(D\times W)}\left\|\nabla_{z_1}f\right\|_{L^2(\Omega_{R,\varepsilon})}\cr
&=\gamma\left\|\nabla_{\bar z_1}f\right\|_{L^2(D\times W)}\,\left\|f\right\|_{L^{4\gamma-2}}^{2\gamma-1}+\gamma\left\|f\right\|_{L^{4\gamma-2}}^{2\gamma-1}\left\|\nabla_{z_1}f\right\|_{L^2(\Omega_{R,\varepsilon})}\cr
&=\gamma\left\|f\right\|_{L^{4\gamma-2}(D\times W)}^{2\gamma-1}\left(\left\|\nabla_{\bar z_1}f\right\|_{L^2(D\times W)}+\left\|\nabla_{z_1}f\right\|_{L^2(D\times W)}\right).\cr
\end{aligned}
$$
We now apply the Cauchy integral formula for nonholomorphic smooth functions to
the function $\Phi=\left(|f|^2\right)^{\gamma}$ to get
$$
\left|\Phi(z_1)\right|\leq\frac{1}{\pi}\int_{\zeta_1\in D}\frac{1}{|\zeta_1-z_1|}\left|\partial_{\bar\zeta_1}\Phi(\zeta_1)\right|
$$
for $z_1\in D$.
We assume that $D\times W$ is contained in $\left({\mathbb D}_{\hat R}\right)^n$ for some $\hat R>1$.
Fix some $0<\eta<1$.  From the potential estimate applied to
$$
h(z,\zeta)=\frac{1}{\pi|z_1-\zeta_1|}
$$
on $D\times W$ with $r=2-\eta$, $p=1$, and $q=2-\eta$,
we have
$$
\begin{aligned}\sup_{x\in\Omega}\|h(x,\cdot)\|_{L^r(D\times W)}
&=\left(\sup_{z\in D\times W}\int_{\zeta\in D\times W}\frac{1}{(\pi|\zeta_1-z_1|)^{2-\eta}}\right)^{\frac{1}{2-\eta}}\cr
&\leq\left((\pi\hat R)^{n-1}\int_{|\zeta_1|<2\hat R}\frac{1}{(\pi|\zeta_1|)^{2-\eta}}\right)^{\frac{1}{2-\eta}}=C_{\hat R},\cr
\end{aligned}
$$
where
$$
C_{\hat R}=\left(\frac{2^{1+\eta}\pi^{n-2+\eta}\hat R^{n-1+\eta}}{\eta}\right)^{\frac{1}{2-\eta}},
$$
and
$$
\begin{aligned}\left\|f\right\|_{L^{2\gamma(2-\eta)}(D\times W)}^{2\gamma}&=
\left\|\left(|f|^2\right)^{\gamma}\right\|_{L^q(D\times W)}\leq C_{\hat R}\int_{D\times W}\left|\partial_{\bar z_1}\left(|f|^2\right)^{\gamma}\right|\cr
&\leq
\gamma\, C_{\hat R}\left\|f\right\|_{L^{4\gamma-2}(D\times W)}^{2\gamma-1}\left(\left\|\nabla_{\bar z_1}f\right\|_{L^2(D\times W)}+\left\|\nabla_{z_1}f\right\|_{L^2(D\times W)}\right).\cr
\end{aligned}
$$
Let $\gamma^*$ be the number which satisfies $4\gamma^*-2=2\gamma(2-\eta)$ so that the above inequality becomes
$$
\left\|f\right\|_{L^{4\gamma^*-2}(D\times W)}^{2\gamma}\leq
\gamma\, C_{\hat R}\left\|f\right\|_{L^{4\gamma-2}(D\times W)}^{2\gamma-1}\left(\left\|\nabla_{\bar z_1}f\right\|_{L^2(D\times W)}+\left\|\nabla_{z_1}f\right\|_{L^2(D\times W)}\right).
$$
This means that
$$
4\gamma^*-2=2\gamma(2-\eta)=4\gamma-2\gamma\eta
$$
and
$$
\gamma^*=\gamma\left(1-\frac{\eta}{2}\right)+\frac{1}{2}.
$$ 
Take a positive integer $N$. We now start with $\gamma=1$ to iterate this procedure until we get to an estimate for $\|f\|_{L^N(D\times W)}$.  We set $\gamma_0=1$ and inductively define
$$
\gamma_{\nu+1}=\gamma_\nu\left(1-\frac{\eta}{2}\right)+\frac{1}{2}
$$
for $0\leq\nu\leq N+2$ so that the relation of $\gamma_\nu$ to $\gamma_{\nu+1}$ is the relation of $\gamma$ to $\gamma^*$.  Then
$$
\gamma_\nu=\frac{1}{\eta}\left(1-(1-\eta)\left(1-\frac{\eta}{2}\right)^\nu\right)
$$
for $0\leq\nu\leq N+2$,
because the induction of going from the $\nu$-th step to the $(\nu+1)$-th step is verified by
$$
\begin{aligned}
\gamma_\nu\left(1-\frac{\eta}{2}\right)+\frac{1}{2}
&=\frac{1}{\eta}\left(\left(1-\frac{\eta}{2}\right)-(1-\eta)\left(1-\frac{\eta}{2}\right)^{\nu+1}\right)+\frac{1}{2}\cr
&=\frac{1}{\eta}\left(1-(1-\eta)\left(1-\frac{\eta}{2}\right)^{\nu+1}\right)-\frac{\ \frac{\eta}{2}\ }{\eta}+\frac{1}{2}\cr
&=\frac{1}{\eta}\left(1-(1-\eta)\left(1-\frac{\eta}{2}\right)^{\nu+1}\right)=\gamma_{\nu+1}.\cr
\end{aligned}
$$
We set $\eta=\frac{2}{N+2}$ and let $\hat\nu$ be the smallest integer which is no less than
$$
\frac{\log 2}{\log\frac{N+2}{N+1}}.
$$
Then
$$
\left(1-\frac{\eta}{2}\right)^{\hat\nu}=\left(\frac{N+1}{N+2}\right)^{\hat\nu}\leq\frac{1}{2}
$$
and
$$
\gamma_{\hat\nu}\geq\frac{1}{\eta}\left(1-\left(1-\frac{\eta}{2}\right)^{\hat\nu}\right)\geq\frac{1}{2\eta}=\frac{N+2}{4}
$$
so that $4\gamma_{\hat\nu}-2\geq N$.  We apply to the last term on the right-hand side of
$$
\left\|f\right\|_{L^{4\gamma_{\nu+1}-2}(D\times W)}^{2\gamma_\nu}\leq
\gamma_\nu C_{\hat R}\left\|f\right\|_{L^{4\gamma_\nu-2}(D\times W)}^{2\gamma_\nu-1}\left(\left\|\nabla_{\bar z_1}f\right\|_{L^2(D\times W)}+\left\|\nabla_{z_1}f\right\|_{L^2(D\times W)}\right)
$$
for $0\leq\nu<\hat\nu$, the arithmetic-geometric-means inequality
$$
AB=\left(\frac{1}{\delta^{\frac{1-\alpha}{\alpha}}}A^{\frac{1}{\alpha}}\right)^\alpha\left(\delta B^{\frac{1}{1-\alpha}}\right)^{1-\alpha}\leq\frac{\alpha}{\delta^{\frac{1-\alpha}{\alpha}}}A^{\frac{1}{\alpha}}+(1-\alpha)\delta B^{\frac{1}{1-\alpha}}
$$
for $0<\alpha<1$ and $\delta, A, B>0$.  We use $\alpha=\frac{2\gamma_\nu-1}{2\gamma_\nu}$ to get
$$
\begin{aligned}\left\|f\right\|_{L^{4\gamma_{\nu+1}-2}(D\times W)}^{2\gamma_\nu}\leq&
\gamma_\nu C_{\hat R}\,\frac{2\gamma_\nu-1}{2\gamma_\nu\delta_\nu^*}\,\left\|f\right\|_{L^{4\gamma_\nu-2}(D\times W)}^{2\gamma_\nu}\cr
&+\gamma_\nu C_{\hat R}\,\frac{\delta_\nu}{2\gamma_\nu}\left(\left\|\nabla_{\bar z_1}f\right\|_{L^2(D\times W)}+\left\|\nabla_{z_1}f\right\|_{L^2(D\times W)}\right)^{2\gamma_\nu}
\end{aligned}
$$
for any $\delta_\nu>0$, where $\delta_\nu^*=\delta_\nu^{\frac{1}{2\gamma_\nu-1}}$.

\bigbreak Since for any nonnegative numbers $A, B, C$ and any $\theta>0$ we have
$$
A^\theta+B^\theta+C^\theta\leq 3\max(A^\theta,B^\theta,C^\theta)=3\left(\max(A, B, C)\right)^\theta\leq 3\left(A+B+C\right)^\theta,
$$
we can write the above inequality as
$$
\begin{aligned}\left\|f\right\|^2_{L^{4\gamma_{\nu+1}-2}(D\times W)}\leq&
3^{{}^{\frac{1}{2\gamma_\nu}}}\gamma_\nu C_{\hat R}\,\frac{2\gamma_\nu-1}{\gamma_\nu\delta_\nu^*}\,\left\|f\right\|^2_{L^{4\gamma_\nu-2}(D\times W)}\cr
&+3^{{}^{\frac{1}{2\gamma_\nu}}}2\gamma_\nu C_{\hat R}\,\frac{\delta_\nu}{\gamma_\nu}\left(\left\|\nabla_{\bar z_1}f\right\|^2_{L^2(D\times W)}+\left\|\nabla_{z_1}f\right\|^2_{L^2(D\times W)}\right)
\end{aligned}
$$
for any $\delta_\nu>0$.  Let $C^\natural$ be the maximum of
$$
3^{{}^{\frac{1}{2\gamma_\nu}}}\gamma_\nu C_{\hat R}\,\frac{2\gamma_\nu-1}{\gamma_\nu\delta_\nu^*}\quad{\rm and}\quad
3^{{}^{\frac{1}{2\gamma_\nu}}}2\gamma_\nu C_{\hat R}\,\frac{\delta_\nu}{\gamma_\nu}
$$
for $0\leq\nu<\hat\nu$ so that we have
$$
\left\|f\right\|^2_{L^{4\gamma_{\nu+1}-2}(D\times W)}\leq
C^\natural\left(\left\|f\right\|^2_{L^{4\gamma_\nu-2}(D\times W)}
+\left\|\nabla_{\bar z_1}f\right\|^2_{L^2(D\times W)}+\left\|\nabla_{z_1}f\right\|^2_{L^2(D\times W)}\right)\leqno{(\natural)_{\varepsilon,\nu}}
$$
Let $C^\sharp=\sum_{j=0}^{\hat\nu-1}(C^\natural)^j$.  Putting together $(\natural)_{\varepsilon,\nu}$ for $0\leq\nu<\hat\nu$, we get
$$
\begin{aligned}\left\|f\right\|^2_{L^{4\gamma_{\hat\nu}-2}(D\times W)}&\leq
(C^\natural)^{\hat\nu}\left\|f\right\|^2_{L^2(D\times W)}
+C^\sharp\left(\left\|\nabla_{\bar z_1}f\right\|^2_{L^2(D\times W)}+\left\|\nabla_{z_1}f\right\|^2_{L^2(D\times W)}\right)\cr
&\leq
(C^\natural)^{\hat\nu}\left\|f\right\|^2_{L^2(D\times W)}
+C^\sharp
\left(2\left\|\bar\nabla_{z_1}f\right\|^2_{L^2(D\times W)}+
\int_{D\times W}
|\Theta_{1\bar 1}|\,|f|^2\right).
\end{aligned}
$$
We now apply H\"older's inequality to get
$$
\left\|f\right\|^2_{L^N(D\times W)}\leq C^\flat\left\|f\right\|^2_{L^{4\gamma_{\hat\nu}-2}(D\times W)}
$$
from $4\gamma_{\hat\nu}-2\geq N$,
where
$$
C^\flat=\left(\pi^n \hat R^{2n}\right)^{\frac{2(4\gamma_{\hat\nu}-2-N)}{4\gamma_{\hat\nu}-2}}.
$$
For
${\mathbf a}=(a_{\alpha\bar\beta})_{\alpha,\beta}$, ${\mathbf b}=(b^\alpha)_\alpha$, and ${\mathbf c}=(c^\beta)_\beta$,
the inequality
$$
\left|\sum_{\alpha,\beta}a_{\alpha\bar\beta\,}b^\alpha\,\overline{c^\beta}\right|
\leq|{\mathbf a}|\,|{\mathbf b}|\,|{\mathbf c}|
$$
follows from taking the square roots of the Cauchy-Schwarz inequality
$$
\left|\sum_{\alpha,\beta}a_{\alpha\bar\beta\,}b^\alpha\,\overline{c^\beta}\right|^2\leq
\left(\sum_{\alpha,\beta}\left|a_{\alpha\bar\beta}\right|^2\right)\left(\sum_{\alpha,\beta}\left|b^\alpha\,\overline{c^\beta}\right|^2\right)
=\left(\sum_{\alpha,\beta}\left|a_{\alpha\bar\beta}\right|^2\right)\left(\sum_\alpha\left|b^\alpha\right|^2\right)\left(\sum_\beta|c^\beta|^2\right),
$$
which, in the case of
$a_{\alpha\bar\beta}=\Theta_{\alpha\bar\beta 1\bar 1}$ and ${\mathbf b}={\mathbf c}=(f^\alpha)_\alpha$, yields
$$
\left|\sum_{\alpha,\beta}\Theta_{\alpha\bar\beta 1\bar 1}f^\alpha\,\overline{f^\beta}\right|\leq|\Theta_{1\bar 1}|\,|f|^2.
$$
Hence, from H\"older's inequality,
$$
\begin{aligned}
&\int_{D\times W}|\Theta_{1\bar 1}|\,|f|^2
\leq\|f\|^2_{L^N(D\times W)}\left(\int_{D\times W}|\Theta|^{\frac{N}{N-2}}\right)^{\frac{N-2}{N}}\cr
&\leq C^\flat\left\|f\right\|^2_{L^{4\gamma_{\hat\nu}-2}(D\times W)}\left(\int_{D\times W}|\Theta|^{\frac{N}{N-2}}\right)^{\frac{N-2}{N}}\cr
&\leq C^\flat\left(\int_{D\times W}|\Theta|^{\frac{N}{N-2}}\right)^{\frac{N-2}{N}}
\left((C^\natural)^{\hat\nu}\left\|f\right\|^2_{L^2(D\times W)}
+C^\sharp
\left(2\left\|\bar\nabla_{z_1}f\right\|^2_{L^2(D\times W)}+
\int_{D\times W}
|\Theta_{1\bar 1}|\,|f|^2\right)\right).\cr
\end{aligned}
$$
Finally, moving, to the left-hand side, the term on the right-hand side which contains the factor
$$
\int_{D\times W}|\Theta_{1\bar 1}|\,|f|^2
$$
and dividing by its coefficients after the movement, we obtain
$$
\begin{aligned}
\int_{D\times W}|\Theta_{1\bar 1}|\,|f|^2
\leq&\left(1-C^\sharp C^\flat\left(\int_{D\times W}|\Theta|^{\frac{N}{N-2}}\right)^{\frac{N-2}{N}}\right)^{-1}\,C^\flat\left(\int_{D\times W}|\Theta_{1\bar 1}|^{\frac{N}{N-2}}\right)\cdot\cr
&\ \ \cdot
\left((C^\natural)^{\hat\nu}\left\|f\right\|^2_{L^2(D\times W)}
+2C^\sharp
\left\|\bar\nabla_{z_1}f\right\|^2_{L^2(D\times W)}\right)
\end{aligned}
$$
if
$$
1-C^\sharp C^\flat\left(\int_{D\times W}|\Theta|^{\frac{N}{N-2}}\right)^{\frac{N-2}{N}}>0.
$$
Suppose
$$
C^\sharp C^\flat\left(\int_{D\times W}|\Theta|^{\frac{N}{N-2}}\right)^{\frac{N-2}{N}}<\frac{1}{2}.
$$
Let
$$
\delta=2\,C^\flat\left(\int_{D\times W}|\Theta|^{\frac{N}{N-2}}\right)^{\frac{N-2}{N}}\left((C^\natural)^{\hat\nu}+2C^\sharp\right).
$$
Then
$$
\int_{D\times W}|\Theta_{1\bar 1}|\,|f|^2\leq\delta\left(\left\|f\right\|^2_{L^2(D\times W)}+
\left\|\bar\nabla_{z_1}f\right\|^2_{L^2(D\times W)}\right)
$$
and
$$
\begin{aligned}
\left\|f\right\|^2_{L^N(D\times W)}&\leq C^\flat\left\|f\right\|^2_{L^{4\gamma_{\hat\nu}-2}(D\times W)}\cr
&\leq
C^\flat(C^\natural)^{\hat\nu}\left\|f\right\|^2_{L^2(D\times W)}\cr
&\quad+C^\flat C^\sharp
\left(2\left\|\bar\nabla_{z_1}f\right\|^2_{L^2(D\times W)}+\delta\left(\left\|f\right\|^2_{L^2(D\times W)}+
\left\|\bar\nabla_{z_1}f\right\|^2_{L^2(D\times W)}\right)\right)\cr
&\leq
C^\flat\left((C^\natural)^{\hat\nu}+\delta C^\sharp\right)\left\|f\right\|^2_{L^2(D\times W)}
+(2+\delta)C^\flat C^\sharp
\left\|\bar\nabla_{z_1}f\right\|^2_{L^2(D\times W)}.\cr
\end{aligned}
$$
Note that the constants $C^\flat$, $C^\#$, $C^\natural$ depend only on $\hat R$, $N$, and $\eta$ and are independent of $D$ and $W$.

\vskip.3in\bigbreak\noindent{\bf\S4}\ \  {\sc Solving $\bar\partial$ Equation to Construct Sections}

\bigbreak\noindent(4.1) {\it Proposition (Solution of Bundle-Valued $\bar\partial$ Equation with $L^2$ Estimates).}  Let $\hat R>1$ and $n\geq 2$ be an integer.  Let $E$ be a holomorphic vector bundle on ${\mathbb D}_{\hat R}^*\times{\mathbb D}_{\hat R}^{n-1}$ with Hermitian metric $h$ such that the curvature tensor $\Theta$ of $h$ satisfies
$$
\int_{{\mathbb D}_{\hat R}^*\times{\mathbb D}_{\hat R}^{n-1}}\left|\Theta\right|^{\frac{N}{N-2}}<\infty
$$
for some integer $N>2$.  Here $|\Theta|$ means the pointwise norm of the curvature tensor $\Theta$ with respect to the metric $h$ and the Euclidean metric of ${\mathbb C}^n$.  For $\gamma>0$ there exists $0<R_\gamma<\hat R$ such that
$$
\int_{{\mathbb D}_{R_\gamma}^*\times{\mathbb D}_R^{n-1}}\left|\Theta\right|^{\frac{N}{N-2}}<\gamma.
$$
We can choose $\gamma$ so small that
$$
C^\sharp C^\flat\gamma^{\frac{N-2}{N}}<\frac{1}{2}.
$$
We will impose more condition on $\gamma$ later.  Arbitrarily choose any $0<\varepsilon<R_\gamma$.  Let $\Omega_{\gamma,\varepsilon}=\left({\mathbb D}_{R_\gamma}-\overline{{\mathbb D}_\varepsilon}\right)\times{\mathbb D}_R^{n-1}$.  Choose any $K>0$ and introduce the additionally twisted Hermitian metric $h_K=e^{-K|z|^2}h$ of $E$.  Let $g=\sum_{j=1}^ng_{\bar j}d\bar z_j$ be a smooth $E$-valued $(0,1)$-form on $\overline{\Omega_{\gamma,\varepsilon}}$ which belongs to the domain of the actual adjoint $\bar\partial^*$ of $\bar\partial$.  The condition that $g$ belongs to the domain of $\bar\partial^*$ means that $g_{\bar 1}$ is identically zero on $\left(\partial\left({\mathbb D}_{R_\gamma}-\overline{{\mathbb D}_\varepsilon}\right)\right)\times{\mathbb D}_R^{n-1}$ and $g_{\bar j}$ is identically zero on
$$
\left({\mathbb D}_{R_\gamma}-\overline{{\mathbb D}_\varepsilon}\right)\times{\mathbb D}_R^{j-2}\times(\partial{\mathbb D}_R)\times{\mathbb D}_R^{n-j}
$$
for $2\leq j\leq n$.  This condition does not depend on which one of the two metrics $h$ and $h_K$ is used to define $\bar\partial^*$.  Since the Levi-form is identically zero on the smooth part of the boundary of $\Omega_{\gamma,\varepsilon}$, it follows from the Bochner-Kodaira formula with boundary term that
$$
\begin{aligned}&\left\|\bar\partial g\right\|_{L^2(h_K,\Omega_{\gamma,\varepsilon})}^2+\left\|\bar\partial_{h_K}^* g\right\|_{L^2(h_K,\Omega_{\gamma,\varepsilon})}^2\cr
&=\left\|\bar\nabla g\right\|_{{L^2(h_K,\Omega_{\gamma,\varepsilon})}}^2+
\int_{\Omega_{\gamma,\varepsilon}}\left<\Theta\,g,\,g\right>\,e^{-K|z|^2}+K\left\|g\right\|^2_{{L^2(h_K,\Omega_{\gamma,\varepsilon})}},\cr
\end{aligned}
$$
where the symbol $h_K$ in both notations $\bar\partial_{h_K}^*$ and $L^2(h_K,\Omega_{\gamma,\varepsilon})$ indicates that the metric $h_K$ is used in their definitions.

\bigbreak We apply the potential estimate derived in (3.2) to $f=g_{\bar 1}$ on $D\times W$ with $D={\mathbb D}_{R_\gamma}-\overline{{\mathbb D}_\varepsilon}$ and $W={\mathbb D}_R^{n-1}$. Then for $2\leq j\leq n$ we apply the potential estimate derived in (3.2) to $f=g_{\bar j}$ for $2\leq j\leq n$ on $W_1\times D\times W_2$ with $D={\mathbb D}_R$ and $W_1=({\mathbb D}_{R_\gamma}-\overline{{\mathbb D}_\varepsilon})\times{\mathbb D}_R^{j-2}$ and $W_2={\mathbb D}_R^{n-j}$, where in the application of the potential estimate, we interchange the roles of $z_1$ and $z_j$ so that $W_1\times D\times W_2$ becomes $D\times W$ with $W=W_1\times W_2$.  We have
$$
\begin{aligned}
\left\|g_{\bar j}\right\|^2_{L^N(\Omega_{\gamma,\varepsilon})}&\leq
C^\flat\left((C^\natural)^{\hat\nu}+\delta C^\sharp\right)\left\|g_{\bar j}\right\|^2_{L^2(\Omega_{\gamma,\varepsilon})}
+(2+\delta)C^\flat C^\sharp
\left\|\bar\nabla_{z_1}g_{\bar j}\right\|^2_{L^2(\Omega_{\gamma,\varepsilon})}\cr
&\leq
C^\flat\left((C^\natural)^{\hat\nu}+\delta C^\sharp\right)\left\|g\right\|^2_{L^2(\Omega_{\gamma,\varepsilon})}
+(2+\delta)C^\flat C^\sharp
\left\|\bar\nabla g\right\|^2_{L^2(\Omega_{\gamma,\varepsilon})},\cr
\end{aligned}
$$
where
$$
\delta=2\,C^\flat\gamma^{\frac{N-2}{N}}\left((C^\natural)^{\hat\nu}+2C^\sharp\right).
$$
For fixed $1\leq j,k\leq n$ we have the following estimate of the curvature term
$$
\begin{aligned}
&\int_{\Omega_{\gamma,\varepsilon}}\left|\sum_{\alpha,\beta}\Theta_{\alpha\bar\beta j\bar k}g^\alpha_{\bar k}\overline{g^\beta_{\bar j}}\right|
\leq \int_{\Omega_{\gamma,\varepsilon}}|\Theta_{j\bar k}|\,|g_{\bar k}|\,|g_{\bar j}|\cr
&\leq\|g_{\bar k}\|_{L^N(\Omega_{\gamma,\varepsilon})}\,\|g_{\bar j}\|_{L^N(\Omega_{\gamma,\varepsilon})}\,\left(\int_{\Omega_{\gamma,\varepsilon}}|\Theta|^{\frac{N}{N-2}}\right)^{\frac{N-2}{N}}\cr
&\leq\gamma^{\frac{N-2}{N}}
\left(C^\flat\left((C^\natural)^{\hat\nu}+\delta C^\sharp\right)\left\|g\right\|^2_{L^2(\Omega_{\gamma,\varepsilon})}+(2+\delta)C^\flat C^\sharp
\left\|\bar\nabla g\right\|^2_{L^2(\Omega_{\gamma,\varepsilon})}\right).\cr
\end{aligned}
$$
Let $\kappa=e^{-K\hat R^2}$.  Then
$$
\begin{aligned}
&\left\|\bar\nabla g\right\|_{{L^2(h_K,\Omega_{\gamma,\varepsilon})}}^2+
\int_{\Omega_{\gamma,\varepsilon}}\left<\Theta\,g,\,g\right>\,e^{-K|z|^2}+K\left\|g\right\|^2_{{L^2(h_K,\Omega_{\gamma,\varepsilon})}}\cr
&\geq\kappa\left\|\bar\nabla g\right\|_{{L^2(h_K,\Omega_{\gamma,\varepsilon})}}^2
-\int_{\Omega_{\gamma,\varepsilon}}\left|\left<\Theta\,g,\,g\right>\right|
+\kappa\, K\left\|g\right\|^2_{{L^2(\Omega_{\gamma,\varepsilon})}}\cr
&\geq\kappa\left\|\bar\nabla g\right\|_{{L^2(h_K,\Omega_{\gamma,\varepsilon})}}^2
+\kappa\, K\left\|g\right\|^2_{{L^2(\Omega_{\gamma,\varepsilon})}}\cr
&\quad-n^2\gamma^{\frac{N-2}{N}}
\left(C^\flat\left((C^\natural)^{\hat\nu}+\delta C^\sharp\right)\left\|g\right\|^2_{L^2(\Omega_{\gamma,\varepsilon})}+(2+\delta)C^\flat C^\sharp
\left\|\bar\nabla g\right\|^2_{L^2(\Omega_{\gamma,\varepsilon})}\right)\cr
&=\left(\kappa-n^2\gamma^{\frac{N-2}{N}}(2+\delta)C^\flat C^\sharp\right)\left\|\bar\nabla g\right\|_{{L^2(h_K,\Omega_{\gamma,\varepsilon})}}^2\cr
&\quad+\left(\kappa\, K-n^2\gamma^{\frac{N-2}{N}}
\left(C^\flat\left((C^\natural)^{\hat\nu}+\delta C^\sharp\right)\right)\right)\left\|g\right\|^2_{{L^2(\Omega_{\gamma,\varepsilon})}}.\cr
\end{aligned}
$$
We impose on $\gamma$ the additional condition that
$$
n^2\gamma^{\frac{N-2}{N}}(2+\delta)C^\flat C^\sharp\leq\frac{\kappa}{2}
$$
and
$$
n^2\gamma^{\frac{N-2}{N}}
\left(C^\flat\left((C^\natural)^{\hat\nu}+\delta C^\sharp\right)\right)\leq\frac{\kappa\,K}{2}
$$
so that
$$
\left\|\bar\partial g\right\|_{L^2(h_K,\Omega_{\gamma,\varepsilon})}^2+\left\|\bar\partial_{h_K}^* g\right\|_{L^2(h_K,\Omega_{\gamma,\varepsilon})}^2
\geq\frac{\kappa\,K}{2}\left\|g\right\|^2_{{L^2(\Omega_{\gamma,\varepsilon})}}\geq\frac{\kappa\,K}{2}\left\|g\right\|_{L^2(h_K,\Omega_{\gamma,\varepsilon})}^2.
$$
By the standard smoothing argument of Friedrichs and H\"ormander and the standard functional analysis argument for solving the $\bar\partial$ equation with $L^2$ estimates, we conclude that for any $E$-valued $\bar\partial$-closed $(0,1)$-form $v$ on $\Omega_{\gamma,\varepsilon}$, there exists a smooth $E$-valued function $u_\varepsilon$ on $\Omega_{\gamma,\varepsilon}$ with
$$
\left\|u_\varepsilon\right\|_{L^2(h_K,\Omega_{\gamma,\varepsilon})}^2\leq\sqrt{\frac{2}{\kappa\,K}}
\left\|v\right\|_{L^2(h_K,\Omega_{\gamma,\varepsilon})}^2.
$$
Since the constants $K,\kappa$ are independent of $\varepsilon$, by the standard argument of passing to limit as $\varepsilon\to 0$, we obtain a solution $u$ of the equation $\bar\partial u=v$ on ${\mathbb D}_{R_\gamma}^*\times{\mathbb D}_R^{n-1}$ with the estimate
$$
\left\|u\right\|_{L^2\left(h_K,{\mathbb D}_{R_\gamma}^*\times{\mathbb D}_R^{n-1}\right)}^2\leq\sqrt{\frac{2}{\kappa\,K}}
\left\|v\right\|_{L^2\left(h_K,{\mathbb D}_{R_\gamma}^*\times{\mathbb D}_R^{n-1}\right)}^2.
$$

\vskip.2in\bigbreak\noindent(4.2) {\it Proposition (Constructing Sections Achieving Prescribed Jet Values).}  Let $q$ be a nonnegative integer.  For any $P_0$ in ${\mathbb D}_{R_\gamma}^*\times{\mathbb D}_R^{n-1}$ and any prescribed jet $F_{P_0}$ at $P_0$ of order $q$, there exists a holomorphic $L^2$ section of $E$ over
${{\mathbb D}_{R_\gamma}^*\times{\mathbb D}_R^{n-1}}$ which achieves the prescribed jet $F_{P_0}$ at $P_0$.

\bigbreak\noindent{\it Proof.}  Let $P_0=(z_1^*,\cdots,z_n^*)$ and choose $0<\varepsilon<|z_1^*|$ such that $R+\varepsilon<\hat R$.  Since $({\mathbb D}_{R_\gamma}-\overline{{\mathbb D}_\varepsilon})\times{\mathbb D}_{R+\varepsilon}^{n-1}$ is Stein, we can find a holomorphic section $s$ of $E$ on $({\mathbb D}_{R_\gamma}-\overline{{\mathbb D}_\varepsilon})\times{\mathbb D}_{R+\varepsilon}^{n-1}$ such that the $q$-jet of $s$ at $P_0$ agrees with the prescribed $q$-jet $F_{P_0}$ at $P_0$.  Let $0\leq\rho(z_1)\leq 1$ be a smooth function on ${\mathbb D}_{R_\gamma}-\overline{{\mathbb D}_\varepsilon}$ with compact support which is identically $1$ on some open neighborhood of $z_1^*$.  Let $v$ be the $E$-valued $(0,1)$-form
$$\frac{(\bar\partial\rho)s}{(z_1-z_1^*)^q}$$
on ${{\mathbb D}_{R_\gamma}^*\times{\mathbb D}_R^{n-1}}$. The $L^2$ norm of $v$ on ${{\mathbb D}_{R_\gamma}^*\times{\mathbb D}_R^{n-1}}$ is finite.  By Proposition(4.1) we can solve the equation $\bar\partial u=v$ on ${{\mathbb D}_{R_\gamma}^*\times{\mathbb D}_R^{n-1}}$ to get a smooth section $u$ of $E$ on ${{\mathbb D}_{R_\gamma}^*\times{\mathbb D}_R^{n-1}}$ with finite $L^2$ norm so that $s-(z_1-z_1^*)^qu$ is a holomorphic $L^2$ section of $E$ over
${{\mathbb D}_{R_\gamma}^*\times{\mathbb D}_R^{n-1}}$ which achieves the prescribed $q$-jet $F_{P_0}$ at $P_0$. Q.E.D.

\vskip.3in\bigbreak\noindent{\bf\S5}\ \  {\sc Metrics with Removable-Singularity Property, Gap-Sheaves and Singularity Subvarieties of Coherent Sheaves}

\bigbreak\noindent(5.1) {\it Two Ingredients Needed to Piece Together Local Extensions of Given Holomorphic Vector Bundle.} In order to prove the Main Theorem, we will start out with a local Thullen-extension setting at a
point of $G$ to construct by Proposition(4.2) (with $q=0$) holomorphic sections to locally extend ${\mathcal O}(E)$ as a reflexive coherent sheaf and then use a sequence of connecting local Thullen-extension settings to get to any point on $Y$.  For such a procedure, we have to handle two difficulties.

\bigbreak One difficulty is that we need to extend the metric to the locally free part of the local extension of ${\mathcal O}(E)$ before we can move on to the next connecting local Thullen-extension setting.  Though a smooth extension of the metric may not be possible, we can still prove that a possibly singular extension of the metric possesses a removable-singularity property which can perform the same task.  We are going to introduce the definition of this removable-singularity property of the metric.

\bigbreak The other difficulty is how to piece together all the local reflexive coherent sheaves obtained from the local extensions of ${\mathcal O}(E)$ to yield a unique well-defined global reflexive coherent sheaf.   We are going to use the gap-sheaves of coherent analytic sheaves to get the global extension.  We will recall from [Siu-Tra] what is needed from the theory of gap-sheaves and the singularity sets of coherent shaves.

\bigbreak\noindent(5.2) {\it Definition of Removable-Singularity Property of Metric.}  Suppose $E$ is a holomorphic vector bundle on ${\mathbb D}^n$ with a smooth Hermitian metric $h$ only defined on  ${\mathbb D}^*\times{\mathbb D}^{n-1}$.  The metric $h$ is said to possess {\it the removable singularity property} at a point $P$ of $\{0\}\times{\mathbb D}^{n-1}$ if any holomorphic section $f$ of $E$ on $U\cap\left({\mathbb D}^*\times{\mathbb D}^{n-1}\right)$ for some open neighborhood $U$ of $P$ in ${\mathbb D}^n$ which is $L^2$ with respect to $h$ can be extended to a holomorphic section germ of $E$ at $P$.

\bigbreak If the Hermitian matrix representing $h$ on $U\cap\left({\mathbb D}^*\times{\mathbb D}^{n-1}\right)$ with respect to some local trivialization of $E$ on $U$ dominates positive smooth metric of $E$ on $U$, then clearly $h$ possesses the removable singularity property at $P$.

\bigbreak In our application, besides the metric $h$ of $E$ possessing the removable-singularity property, the metric $h^{-1}$ for the dual $E^*$ of $E$ induced by $h$ needs also to possess the removable singularity property.

\bigbreak\noindent(5.3) {\it Gap-Sheaves and Singularity Sets of Coherent Sheaves.}  Suppose ${\mathcal F}$ is a coherent sheaf on a domain $\Omega$ in ${\mathbb C}^n$.  For $0\leq m\leq n$ the singularity subvariety $S_m({\mathcal F})$ is defined as the set of points $x$ where the homological codimension ${\rm codh}_x{\mathcal F}$ of ${\mathcal F}$ is $\leq m$ (see [Siu-Tra, p.26]).  The sheaf ${\mathcal F}$ is locally free at $x$ if and only if ${\rm codh}_x{\mathcal F}=n$.  By Theorem(1.11) on p.31 of [Siu-Tra], the complex dimension of the subvariety $S_m({\mathcal F})$ is always $\leq m$.  The {\it absolute gap-sheaf} ${\mathcal R}^0_d{\mathcal F}$ (with local subvarieties of complex dimension $\leq d$ as gaps) is defined by the presheaf which assigns to an open subset $Q$ of $\Omega$ the vector space ${\mathcal P}_Q$ with the following definition.
An element $s$ of ${\mathcal P}_Q$ is a holomorphic section of ${\mathcal F}$ over $Q-Z_s$ for some (possibly empty) complex subvariety $Z_s$ of complex dimension $\leq d$ in $Q$ (see [Siu-Tra, p.76]).  The sheaf ${\mathcal R}^0_d{\mathcal F}$ is coherent if and only if $\dim S_{d+1}({\mathcal F})\leq d$ (see Corollary(3.9) on p.79 of [Siu-Tra]).  For a reflexive ${\mathcal F}$, the gap-sheaf ${\mathcal R}^0_{n-2}{\mathcal F}$ is always equal to ${\mathcal F}$ and as a consequence $\dim S_{n-1}({\mathcal F})\leq n-2$, which means that ${\mathcal F}$ is locally free outside a subvariety of codimension $\geq 2$ in $\Omega$.

\bigbreak\noindent(5.4) {\it Proposition (Reflexive Sheaf Constructed from Extension of Locally Free Sheaves).}  Let $P^*$ a point of $\{0\}\times{\mathbb D}^{n-1}$ and $U$ be a polydisk open neighborhood of $P^*$ in ${\mathbb D}^n$.
Suppose $E$ is a holomorphic vector bundle on $\left({\mathbb D}^*\times{\mathbb D}^{n-1}\right)\cup U$ and $h$ is a smooth Hermitian metric of $E$ on ${\mathbb D}^*\times{\mathbb D}^{n-1}$ (not assumed to be defined on $U$) such that
\begin{itemize}
\item[(a)] the pointwise norm of the curvature tensor $\Theta_h$ with respect to $h$, as a nonnegative function on ${\mathbb D}^*\times{\mathbb D}^{n-1}$, is locally $L^p$ on ${\mathbb D}^n$ for some $p>1$, and
\item[(b)] both the metric $h$ of $E$ and the metric $h^{-1}$ of $E^*$ possess the removable singularity property at some point $P^*$ of $(\{0\}\times{\mathbb D}^{n-1})\cap U$.
\end{itemize}
Then ${\mathcal O}(E)$ can be uniquely extended to a reflexive coherent sheaf ${\mathcal F}$ on ${\mathbb D}^n$.  The unique extension sheaf ${\mathcal F}$ is defined
by the presheaf which assigns to an open subset $Q$ of ${\mathbb D}^n$ the vector space ${\mathcal P}_Q$ with the following definition.  An element $\sigma$ of ${\mathcal P}_Q$ is a holomorphic section of $E$ over $Q\cap\{z_1\not=0\}$ such that the pointwise norm $|\sigma|$ of $\sigma$ with respect to $h$ as a function on $Q\cap\{z_1\not=0\}$ is locally $L^2$ on $Q-Z_s$ for some (possibly empty) complex subvariety $Z_\sigma$ of complex dimension $\leq n-2$ in $Q\cap\{z_1=0\}$.  Moreover, the metrics $h$ and $h^{-1}$ possess the removable singularity property at every point of $(\{0\}\times{\mathbb D}^{n-1})-S_{n-1}({\mathcal F})$.

\bigbreak\noindent{\it Proof.} There exists $0<\gamma<1$ such that for any $0<\varepsilon<\gamma$ and any $q\in{\mathbb N}$, any prescribed $q$-jet of the fiber of $E$ at any prescribed point can be achieved by some global $L^2$ holomorphic sections of $E$ on ${\mathbb D}_\gamma^*\times{\mathbb D}_{1-\varepsilon}^{n-1}$.  By Condition (b) and Hartogs's extension theorem for holomorphic functions, global $L^2$ holomorphic sections of $E$ on ${\mathbb D}_\gamma^*\times{\mathbb D}_{1-\varepsilon}^{n-1}$ can be extended to holomorphic sections of ${\mathcal O}(E)$ on $\left({\mathbb D}_\gamma^*\times{\mathbb D}_{1-\varepsilon}^{n-1}\right)\cup U$.

\bigbreak Choose a point $P^\flat$ of ${\mathbb D}_\gamma^*\times{\mathbb D}_{1-\varepsilon}^{n-1}$.  We can apply the argument to the dual bundle $E^*$ to get global $L^2$ holomorphic sections $t_1,\cdots,t_\ell$ of $E^*$ over
$\left({\mathbb D}_\gamma^*\times{\mathbb D}_{1-\varepsilon}^{n-1}\right)\cup U$ whose restrictions to ${\mathbb D}_\gamma^*\times{\mathbb D}_{1-\varepsilon}^{n-1}$ have finite $L^2$ norm and which generate the stalk of $E$ at $P^\flat$.
We have a sheaf-homomorphism $\phi:{\mathcal O}(E)\to{\mathcal O}^\ell$ on $\left({\mathbb D}_\gamma^*\times{\mathbb D}_{1-\varepsilon}^{n-1}\right)\cup U$ defined by $t_1,\cdots,t_\ell$, which is a sheaf-monomorphism on some open neighborhood of $P^\flat$ and hence a sheaf-monomorphism on all of $\left({\mathbb D}_\gamma^*\times{\mathbb D}_{1-\varepsilon}^{n-1}\right)\cup U$, because the support of the kernel of $\phi$ is a complex subvariety of $\left({\mathbb D}_\gamma^*\times{\mathbb D}_{1-\varepsilon}^{n-1}\right)\cup U$ not containing $P^\flat$ and is therefore empty on account of the local freeness of ${\mathcal O}(E)$.

\bigbreak Let $S$ be the set of all global $L^2$ holomorphic sections $s$ of $E$ on
$\left({\mathbb D}_\gamma^*\times{\mathbb D}_{1-\varepsilon}^{n-1}\right)\cup U$ whose restrictions to ${\mathbb D}_\gamma^*\times{\mathbb D}_{1-\varepsilon}^{n-1}$ have finite $L^2$ norm.  By Hartogs's extension theorem for holomorphic functions, the image $\phi(s)=\left(\left<t_1,\,s\right>,\cdots,\left<t_\ell,\,s\right>\right)$ of an $s$ element of $S$ can be extended to an $\ell$-tuple $\hat s$ of holomorphic functions on ${\mathbb D}_\gamma\times{\mathbb D}_{1-\varepsilon}^{n-1}$.  Let ${\mathcal L}$ be the subsheaf of ${\mathcal O}^\ell$ on ${\mathbb D}_\gamma\times{\mathbb D}_{1-\varepsilon}^{n-1}$ which is generated by $\hat s$ for $s\in S$.  Let ${\mathcal F}$ be the double dual $${\mathcal Hom}_{\mathcal O}\left({\mathcal Hom}_{\mathcal O}\left({\mathcal L},\,{\mathcal O}\right),\,{\mathcal O}\right)$$
of ${\mathcal L}$. The reflexive sheaf ${\mathcal F}$ on ${\mathbb D}_\gamma\times{\mathbb D}_{1-\varepsilon}^{n-1}$ extends the locally free sheaf ${\mathcal O}(E)$ on $\left({\mathbb D}_\gamma^*\times{\mathbb D}_{1-\varepsilon}^{n-1}\right)\cup U$.

\bigbreak Take a point $P^\#$ of $(\{0\}\times{\mathbb D}^{n-1})-S_{n-1}({\mathcal F})$.  We have to show that both $h$ and $h^{-1}$ possess the removable singularity property at $P^\#$.  Let $U^\#$ be an open polydisk centered at $P^\#$ of polyradii $(r_1,\cdots,r_n)$ with each $r_j<1$ for $1\leq j\leq n$ such that ${\mathcal F}$ is locally free on $U^\#$.   Suppose $f$ is a $L^2$ holomorphic section of ${\mathcal F}$ on $U^\#\cap\{z_1\not=0\}$.  We want to show that $f$ is a holomorphic section of ${\mathcal F}$ on $U^\#$.

\bigbreak We choose a new affine coordinate system $w$ such that
\begin{itemize}
\item[(i)] $w_1=z_1$,
\item[(ii)] $w_j$ vanishes at both $P^*$ and $P^\#$ for $3\leq j\leq n$,
\item[(iii)] $|w_j|<1$ at both $P^*$ and $P^\#$ for $2\leq j\leq n$,
\item[(iii)] $\{|w_2|\leq 1,\cdots,|w_n|\leq 1\}$ is contained in $\{|z_2|<1,\cdots,|z_n|<1\}$,
\item[(iv)] some polydisk in $w$ centered at $P^\#$ of polyradii $(r_1^\prime,\cdots,r_n^\prime)$ with $r_1^\prime=r_1$ and $r_j<1$ for $2\leq j\leq n$ is contained in $U^\#$.
\end{itemize}
Let $0\leq\rho(w_2)\leq 1$ be a smooth function on $\{|w_2-w_2(P^\#)|<r^\prime_2\}$ with compact support which is identically $1$ on $\{|w_2-w_2(P^\#)|\leq\frac{r^\prime_2}{2}\}$.  Take an arbitrary point $P^\flat$ with $w_1(P^\flat)=0$ and $|w_2(P^\flat)-w_2(P^\#)|\leq\frac{r^\prime_2}{2}$ and $|w_j(P^\flat)-w_j(P^\#)|\leq r_j^\prime$ for $3\leq j\leq n$.  By Proposition(4.1), for some $0<\gamma^\prime<r_1$ and any $\varepsilon>0$ with
$|w_2(P^*)|<1-\varepsilon$ and $|w_2(P^\flat)|<1-\varepsilon$, the equation
$$
\bar\partial u=\frac{(\bar\partial\rho)f}{w_2-w_2(P^\flat)}
$$
on $$\Omega^\prime:=\{0<|w_1|<\gamma^\prime,|w_2|<1-\varepsilon,|w_3|<r_3^\prime\cdots,|w_n|<r_n^\prime\}$$ can be solved for an $L^2$ smooth section $u$ of $E$ on $\Omega^\prime$.  Then $f-(w_2-w_2(P^\flat))u$ is a holomorphic $L^2$ section of $E$ on $\Omega^\prime$.  By Condition (b), $f-(w_2-w_2(P^\flat))u$ is holomorphic on an open neighborhood of $P^*$ in ${\mathbb C}^n$.  By Hartogs's extension theorem for holomorphic functions, we conclude that $f-(w_2-w_2(P^\flat))u$ is a holomorphic section of the locally free sheaf ${\mathcal F}$ on $\Omega^\prime-S_{n-1}({\mathcal F})$.
In particular, the restriction of $f$ to $\{w_2=w_2(P^\flat)\}\cap U^\#$ is holomorphic.  Since $P^\flat$ is arbitrarily chosen under the condition $w_1(P^\flat)=0$ and $|w_2(P^\flat)-w_2(P^\#)|\leq\frac{r^\prime_2}{2}$ and $|w_j(P^\flat)-w_j(P^\#)|\leq r_j^\prime$ for $3\leq j\leq n$, it follows that $f$ is holomorphic on some open neighborhood of $P^\#$ in ${\mathbb C}^n$.  This shows that $h$ possesses the removable singularity property at $P^\#$.  By applying the argument to $E^*$ instead of $E$, we conclude that $h^{-1}$ also possesses the removable singularity property at $P^\#$.

\bigbreak The removable singularity property of $h$ at any point $P^\flat$ of $(\{0\}\times{\mathbb D}^{n-1})-S_{n-1}({\mathcal F})$, together with the generation of the subsheaf ${\mathcal L}$ by $\phi(s)$ for $s\in S$, implies that ${\mathcal F}$ is defined by the presheaf which assigns to an open subset $Q$ of $X$ the vector space ${\mathcal P}_Q$ with the following definition. An element $s$ of ${\mathcal P}_Q$ is a holomorphic section of $Q\cap\{z_1\not=0\}$ such that for some (possibly empty) complex subvariety $Z_s$ of complex dimension $\leq n-2$ in $Q\cap\{z_1=0\}$, the pointwise norm $|s|$ of $s$ as a function of $Q\cap\{z_1\not=0\}$ is locally $L^2$ on $Q$ at points of $(Q\cap\{z_1=0\})-Z_s$.
Q.E.D.

\bigbreak We use Proposition(5.4) to yield the following Theorem(5.5) which is slightly more general than the Main Theorem.  The reason for the slightly more general formulation is to facilitate the proof of the Main Theorem.  For every point $P$ of $Y$ we choose a local coordinate system $z_1,\cdots,z_n$ centered at $P$ such that the intersection of $Y$ with the unit polydisk $U_P:={\mathbb D}^n$ of the coordinate system is $\{z_1=0\}$ in ${\mathbb D}^n$.  For any given regular point $P$ of $Y$, we can start out with some regular point $P_0$ of $G\cap Y$ and construct points $P_1,\cdots,P_m$ of $Y$ with $P_m=P$ such that $P_{j+1}$ belongs to $U_{P_j}$ for $1\leq j<m$.  Proposition(5.4) implies by induction on $1\leq j\leq m$ that both $h$ and $h^{-1}$ possess the removable singularity property at $P_j$ and the reflexive coherent sheaf ${\mathcal F}|_{U_{P_j}}$ extends the sheaf ${\mathcal O}(E)|_{U_{P_j}}$.

\bigbreak\noindent(5.5) {\it Theorem (Extension of Bundles with Metrics Possessing Removable-Singularity Property).}  Let $X$ be complex manifold of complex dimension $n\geq 2$ and $Y$ be a connected nonsingular complex hypersurface in $X$ and $G$ be an open subset of $X$ which intersects $Y$.  Let $E$ be a holomorphic vector bundle on $(X-Y)\cup G$ and $h$ be a smooth Hermitian metric of $E|_{X-Y}$ (not assumed to be defined on $G$) such that the pointwise norm of the curvature $\Theta_h$ of $h$, with respect to $h$ and any smooth Hermitian metric of $X$, is locally $L^p$ for some $p>1$ as a function on $X$.  Assume that both $h$ and $h^{-1}$ possess the removable singularity property at every point of $G\cap Y$.  Then $E$ can be extended to a coherent sheaf ${\mathcal F}$ on $X$, which is reflexive and therefore unique.  Moreover, ${\mathcal F}$ can be defined by the presheaf
which assigns to an open subset $Q$ of $X$ the vector space ${\mathcal P}_Q$ with the following definition.  An element $s$ of ${\mathcal P}_Q$ is a holomorphic section of $Q-Y$ such that for some (possibly empty) complex subvariety $Z_s$ of complex dimension $\leq n-2$ in $Q\cap Y$, the pointwise norm $|s|$ of $s$ as a function of $Q-Y$ is locally $L^2$ on $Q-Z_s$.

\vskip.4in\bigbreak\noindent{\sc References}

\bigbreak\noindent [Ang-Siu] Urban Angehrn and Yum-Tong Siu,
Effective freeness and point separation for adjoint bundles.
{\it Invent. Math.} \textbf{122} (1995), 291 -- 308.

\medbreak\noindent[Ban] Shigetoshi Bando,
Removable singularities for holomorphic vector bundles.
{\it Tohoku Math. J.} \textbf{43} (1991), 61 -- 67.

\medbreak\noindent[Ban-Siu] Shigetoshi Bando and Yum-Tong Siu,
Stable sheaves and Einstein-Hermitian metrics. {\it Geometry and analysis on complex manifolds}, 39 -- 50, World Sci. Publ., River Edge, NJ, 1994.

\medbreak\noindent[Bis] Errett Bishop,
Conditions for the analyticity of certain sets.
{\it Michigan Math. J.} \textbf{11} (1964), 289 -- 304.

\medbreak\noindent[Dem1] Jean-Pierre Demailly,
A numerical criterion for very ample line bundles.
{\it J. Differential Geom.} \textbf{37} (1993), 323 -- 374.

\medbreak\noindent[Dem2] Jean-Pierre Demailly,
Effective bounds for very ample line bundles.
{\it Invent. Math.} \textbf{124} (1996), 243 -- 261.

\medbreak\noindent[Fis] Gerd Fischer,
Lineare Faserr\"aume und koh\"arente Modulgarben \"uber komplexen R\"aumen.
{\it Arch. Math. (Basel)} \textbf{18} (1967), 609 -- 617.

\medbreak\noindent[Fuj] Takao Fujita, On polarized manifolds whose adjoint bundles are not semipositive. In Algebraic Geometry, Sendai, {\it Advanced Studies in Pure Math.}
\textbf{10} (1987), 167 -- 178.

\medbreak\noindent[Gil-Tru] David Gilbarg and Neil S. Trudinger,
{\it Elliptic Partial Differential Equations of Second Order}.
Second edition. Springer-Verlag, Berlin, 1983.

\medbreak\noindent[Gra] Hans Grauert, \"Uber Modifikationen und exzeptionelle analytische Mengen. {\it Math. Ann.} \textbf{146},
331 -- 368 (1962).

\medbreak\noindent[Gro] Alexander Grothendieck, Techniques de construction en g\'eom\'etrie analytique. {\it S\'eminaire Henri
Cartan}, \textbf{13c} ann\'ee 1960/61).

\medbreak\noindent[Harv] Reese Harvey,
Removable singularities for positive currents.
{\it Amer. J. Math.} \textbf{96} (1974), 67 -- 78.

\medbreak\noindent[Hart] Friedrich Hartogs,
\"Uber die aus den singul\"aren Stellen einer analytischen Funktion mehrerer Ver\"anderlichen bestehenden Gebilde.
{\it Acta Math.} \textbf{32} (1909), 57 -- 79.

\medbreak\noindent[Hei] Gordon Heier,
Effective freeness of adjoint line bundles.
{\it Doc. Math.} \textbf{7} (2002), 31 -- 42.

\medbreak\noindent[Koh] Joseph J. Kohn,
Subellipticity of the $\bar\partial$-Neumann problem on pseudo-convex domains: sufficient conditions.
{\it Acta Math.} \textbf{142} (1979), 79 -- 122.

\medbreak\noindent[Nad] Alan Nadel,
Multiplier ideal sheaves and K\"ahler-Einstein metrics of positive scalar curvature.
{\it Ann. of Math.} \textbf{132} (1990), 549 -- 596.

\medbreak\noindent[Shi1] Bernard Shiffman,
Extension of positive line bundles and meromorphic maps.
{\it Invent. Math.} \textbf{15} (1972), 332 -- 347.

\medbreak\noindent[Shi2] Bernard Shiffman,
Extension of positive holomorphic line bundles.
{\it Bull. Amer. Math. Soc.} \textbf{77} (1971), 1091 -- 1093.

\medbreak\noindent[Sib] Nessim Sibony,
Quelques probl\`mes de prolongement de courants en analyse complexe.
{\it Duke Math. J.} \textbf{52} (1985), 157 -- 197.

\medbreak\noindent[Siu1] Yum-Tong Siu,
The complex-analyticity of harmonic maps and the strong rigidity of compact K\"ahler manifolds.
{\it Ann. of Math.} \textbf{112} (1980), 73 -- 111.

\medbreak\noindent[Siu2] Yum-Tong Siu,
Some recent results in complex manifold theory related to vanishing theorems for the semipositive case. Workshop Bonn 1984 (Bonn, 1984), 169 -- 192,
{\it Lecture Notes in Math.}, \textbf{1111}, Springer, Berlin, 1985.  (see pp.184-187)

\medbreak\noindent[Siu3]  Yum-Tong Siu,
A Hartogs type extension theorem for coherent analytic sheaves.
{\it Ann. of Math.} \textbf{93} (1971), 166 -- 188.

\medbreak\noindent[Siu4]  Yum-Tong Siu,
A Thullen type theorem on coherent analytic sheaf extension.
{\it Rice Univ. Stud.} \textbf{56} (1970), 187 -- 197 (1971).

\medbreak\noindent[Siu5]  Yum-Tong Siu,
Effective very ampleness.
{\it Invent. Math.} \textbf{124} (1996), 563 -- 571.

\medbreak\noindent[Siu-Tra] Yum-Tong Siu and G\"unther Trautmann,
Gap-sheaves and extension of coherent analytic subsheaves.
{\it Lecture Notes in Mathematics}, Vol. \textbf{172} Springer-Verlag, Berlin-New York 1971.

\medbreak\noindent[Sko] Henri Skoda,
Prolongement des courants, positifs, ferm\'es de masse finie.
{\it Invent. Math.} \textbf{66} (1982), 361 -- 376.

\medbreak\noindent[Thu] Peter Thullen,
\"Uber die wesentlichen Singularit\"aten analytischer Funktionen und Fl\"achen im Raume von n komplexen Ver\"anderlichen.
{\it Math. Ann.} \textbf{111} (1935), 137 -- 157.

\medbreak\noindent[Tra1] G\"unther Trautmann,
Fortsetzung lokal-freier Garben \"uber i-dimensionale Singularit\"atenmengen.
{\t Ann. Scuola Norm. Sup. Pisa Cl. Sci.} \textbf{23} (1969), 155 -- 184.

\medbreak\noindent[Tra2] G\"unther Trautmann,
Ein Kontinuit\"atssatz f\"ur die Fortsetzung koh\"arenter analytischer Garben.
{\it Arch. Math. (Basel)} \textbf{18} (1967), 188 -- 196.

\vskip.2in\bigbreak\noindent{\it Author's affiliation:} Department of Mathematics, Harvard University, Cambridge, MA 02493, USA.

\medbreak\noindent{\it Author's email:} siu@math.harvard.edu

\end{document}